\documentclass{amsart}
\usepackage{amsmath}
\usepackage{amsfonts}

\newtheorem{theorem}{Theorem}
\theoremstyle{plain}

\newtheorem{corollary}{Corollary}

\newtheorem{lemma}{Lemma}

\numberwithin{equation}{section}
\input{tcilatex}

\begin{document}
\title[\textbf{Frobenius-Euler numbers and polynomials}]{\textbf{A note on
the Frobenius-Euler Numbers and polynomials associated with Bernstein
polynomials}}
\author[\textbf{S. Araci}]{\textbf{Serkan Araci}}
\address{\textbf{University of Gaziantep, Faculty of Science and Arts,
Department of Mathematics, 27310 Gaziantep, TURKEY}}
\email{\textbf{mtsrkn@hotmail.com}}
\author[\textbf{M. Acikgoz}]{\textbf{Mehmet Acikgoz}}
\address{\textbf{University of Gaziantep, Faculty of Science and Arts,
Department of Mathematics, 27310 Gaziantep, TURKEY}}
\email{\textbf{acikgoz@gantep.edu.tr}}
\subjclass[2000]{ 05A10, 11B65, 28B99, 11B68, 11B73.}
\keywords{Frobenius-Euler numbers and polynomials, Bernstein polynomials,
fermionic $p$-adic integral on $%
\mathbb{Z}
_{p}$.}

\begin{abstract}
The present paper deals with Bernstein polynomials and Frobenius-Euler
numbers and polynomials. We apply the method of generating function and
fermionic $p$-adic integral representation on $%
\mathbb{Z}
_{p}$, which are exploited to derive further classes of Bernstein
polynomials and Frobenius-Euler numbers and polynomials. To be more precise
we summarize our results as follows, we obtain some combinatorial relations
between Frobenius-Euler numbers and polynomials. Furthermore, we derive an
integral representation of Bernstein polynomials of degree $n$ on $%
\mathbb{Z}
_{p}$. Also we deduce a fermionic $p$-adic integral representation of
product Bernstein polynomials of different degrees $n_{1},n_{2},\cdots $\ on 
$%
\mathbb{Z}
_{p}$ and show that it can be written with Frobenius-Euler numbers which
yields a deeper insight into the effectiveness of this type of
generalizations. Our applications possess a number of interesting properties
which we state in this paper
\end{abstract}

\maketitle

\section{\textbf{Introduction and Notations}}

Let $p$ be a fixed odd prime number. Throughout this paper we use the
following notations. By $%
\mathbb{Z}
_{p}$ we denote the ring of $p$-adic rational integers, $%
\mathbb{Q}
$ denotes the field of rational numbers, $%
\mathbb{Q}
_{p}$ denotes the field of $p$-adic rational numbers, and $%
\mathbb{C}
_{p}$ denotes the completion of algebraic closure of $%
\mathbb{Q}
_{p}$. Let $%
\mathbb{N}
$ be the set of natural numbers and $%
\mathbb{N}
^{\ast }=%
\mathbb{N}
\cup \left\{ 0\right\} $. The $p$-adic absolute value is defined by 
\begin{equation*}
\left\vert p\right\vert _{p}=\frac{1}{p}\text{.}
\end{equation*}
In this paper, we assume $\left\vert q-1\right\vert _{p}<1$ as an
indeterminate. In [17-19], let $UD\left( 
\mathbb{Z}
_{p}\right) $ be the space of uniformly differentiable functions on $%
\mathbb{Z}
_{p}$. For $f\in UD\left( 
\mathbb{Z}
_{p}\right) $, the fermionic $p$-adic integral on $%
\mathbb{Z}
_{p}$ is defined by T. Kim:%
\begin{equation}
I_{-1}\left( f\right) =\int_{%
\mathbb{Z}
_{p}}f\left( \xi \right) d\mu _{-1}\left( \xi \right) =\lim_{N\rightarrow
\infty }\sum_{\xi =0}^{p^{N}-1}f\left( \xi \right) \left( -1\right) ^{\xi }%
\text{.}  \label{equation 1}
\end{equation}

From (\ref{equation 1}), we have well known the following equality:%
\begin{equation}
I_{-1}\left( f_{1}\right) +I_{-1}\left( f\right) =2f\left( 0\right)
\label{equation 4}
\end{equation}

here $f_{1}\left( x\right) :=f\left( x+1\right) $ (for details, see[3-24]).

Let $C\left( \left[ 0,1\right] \right) $ be the space of continuous
functions on $\left[ 0,1\right] $. For $C\left( \left[ 0,1\right] \right) $,
the Bernstein operator for $f$ is defined by 
\begin{equation*}
B_{n}\left( f,x\right) =\sum_{k=0}^{n}f\left( \frac{k}{n}\right)
B_{k,n}\left( x\right) =\sum_{k=0}^{n}f\left( \frac{k}{n}\right) \binom{n}{k}%
x^{k}\left( 1-x\right) ^{n-k}
\end{equation*}

where $n,$ $k\in 
\mathbb{Z}
_{+}:=\left\{ 0,1,2,3,...\right\} $. Here $B_{k,n}\left( x\right) $ is
called Bernstein polynomials, which are defined by 
\begin{equation}
B_{k,n}\left( x\right) =\binom{n}{k}x^{k}\left( 1-x\right) ^{n-k}\text{, }%
x\in \left[ 0,1\right]  \label{equation 2}
\end{equation}

(for more informations on this subject, see [1-6, 11, 14, 15, 17, 21-24])

In \cite{J. Choi}, as is well known, Frobenius-Euler polynomials are defined
by means of the following generating function:%
\begin{equation}
\sum_{n=0}^{\infty }H_{n}\left( u,x\right) \frac{t^{n}}{n!}=e^{H\left(
u,x\right) t}=\frac{1-u}{e^{t}-u}e^{xt}\text{.}  \label{equation 19}
\end{equation}

where the usual convention about replacing $H^{n}\left( u,x\right) $ by $%
H_{n}\left( u,x\right) $. For $x=0$ in (\ref{equation 19}), we have to $%
H_{n}\left( u,0\right) :=H_{n}\left( u\right) $, which is called
Frobenius-Euler numbers. Then, we can write the following%
\begin{equation}
e^{H\left( u\right) t}=\sum_{n=0}^{\infty }H_{n}\left( u\right) \frac{t^{n}}{%
n!}=\frac{1-u}{e^{t}-u}\text{.}  \label{equation 20}
\end{equation}

By (\ref{equation 19}) and (\ref{equation 20}), we easily see the following
applications:%
\begin{eqnarray*}
e^{\left( H\left( u\right) +1\right) t}-ue^{H\left( u\right) t} &=&1-u \\
\sum_{n=0}^{\infty }\left[ \left( H\left( u\right) +1\right)
^{n}-uH_{n}\left( u\right) \right] \frac{t^{n}}{n!} &=&1-u
\end{eqnarray*}

After these applications, we derive the following Lemma.

\begin{lemma}
For $\left\vert u\right\vert >1$ and $n\in 
\mathbb{Z}
_{+}:=%
\mathbb{N}
\tbigcup \left\{ 0\right\} $, we have%
\begin{equation}
\left( H\left( u\right) +1\right) ^{n}-uH_{n}\left( u\right) =\left\{ \QATOP{%
1-u\text{, if }n=0}{0\text{, \ \ \ \ \ \ if }n\neq 0.}\right.
\label{equation 21}
\end{equation}
\end{lemma}

In this paper, we obtained some relations between the Frobenius-Euler
numbers and polynomials and the Bernstein polynomials. From these relations,
we derive some interesting identities on the Frobenius-Euler numbers.

\section{\textbf{On the Frobenius-Euler numbers and polynomials}}

Let us take $f\left( x\right) =u^{x}e^{tx}$ in (\ref{equation 1}), by (\ref%
{equation 4}), we see that 
\begin{equation}
\int_{%
\mathbb{Z}
_{p}}u^{\eta }e^{\eta t}d\mu _{-1}\left( \eta \right) =\frac{2}{1+u}%
H_{n}\left( -u^{-1}\right) \text{.}  \label{equation 22}
\end{equation}

By (\ref{equation 19}) and (\ref{equation 22}), we have the following
theorem.

\begin{theorem}
\begin{equation}
\int_{%
\mathbb{Z}
_{p}}u^{\eta }\left( x+\eta \right) ^{n}d\mu _{-1}\left( \eta \right) =\frac{%
2}{u+1}H_{n}\left( -u^{-1},x\right) \text{.}  \label{equation 23}
\end{equation}
\end{theorem}

By applying some combinatorial techniques in (\ref{equation 23}), we derive
the following%
\begin{equation*}
\int_{%
\mathbb{Z}
_{p}}u^{\eta }\left( x+\eta \right) ^{n}d\mu _{-1}\left( \eta \right)
=\sum_{k=0}^{n}\binom{n}{k}x^{n-k}\left\{ \int_{%
\mathbb{Z}
_{p}}u^{\eta }\eta ^{k}d\mu _{-1}\left( \eta \right) \right\} \text{.}
\end{equation*}

So, from above, we have the well known identity%
\begin{equation}
H_{n}\left( -u^{-1},x\right) =\sum_{k=0}^{n}\binom{n}{k}x^{n-k}H_{k}\left(
-u^{-1}\right) =\left( H\left( -u^{-1}\right) +x\right) ^{n}\text{.}
\label{equation 24}
\end{equation}

by using the $umbral$(symbolic) convention $H^{n}\left( u\right)
:=H_{n}\left( u\right) $.

The Frobenius-Euler polynomials have to symmetric properties, which is shown
by Choi $et$ $al$. in \cite{J. Choi}, as follows: 
\begin{equation*}
H_{n}\left( -u^{-1},1-x\right) =\left( -1\right) ^{n}H_{n}\left(
-u^{-1},x\right) \text{.}
\end{equation*}%
\qquad For $n\in 
\mathbb{N}
$, by (\ref{equation 24}), Choi $et$ $al$. derived the following equality:%
\begin{equation}
u^{2}H_{n}\left( -u^{-1},2\right) =u^{2}+u+H_{n}\left( -u^{-1}\right) \text{.%
}  \label{equation 6}
\end{equation}

From (\ref{equation 23}) and (\ref{equation 6}), we easily see that%
\begin{eqnarray}
&&\int_{%
\mathbb{Z}
_{p}}u^{\eta }\left( 1-\eta \right) ^{n}d\mu _{-1}\left( \eta \right)
\label{equation 13} \\
&=&\left( -1\right) ^{n}\int_{%
\mathbb{Z}
_{p}}u^{\eta }\left( \eta -1\right) ^{n}d\mu _{-1}\left( \eta \right)  \notag
\\
&=&\frac{2}{u+1}\left( -1\right) ^{n}H_{n}\left( -u^{-1},-1\right)  \notag \\
&=&\frac{2}{u+1}H_{n}\left( -u^{-1},2\right) .  \notag
\end{eqnarray}

Thus, we obtain the following Theorem.

\begin{theorem}
The following identity 
\begin{equation}
\int_{%
\mathbb{Z}
_{p}}u^{\eta }\left( 1-\eta \right) ^{n}d\mu _{-1}\left( \eta \right) =\frac{%
2}{u+1}H_{n}\left( -u^{-1},2\right)  \label{equation 15}
\end{equation}%
is true.
\end{theorem}

Let $n\in 
\mathbb{N}
$. By expression of (\ref{equation 6}) and (\ref{equation 15}), we get%
\begin{equation}
\int_{%
\mathbb{Z}
_{p}}u^{\eta }\left( 1-\eta \right) ^{n}d\mu _{-1}\left( \eta \right) =\frac{%
2}{u+1}+\frac{2}{u^{2}+u}+\frac{2}{u^{3}+u}H_{n}\left( -u^{-1}\right) .
\label{equation 8}
\end{equation}

From (\ref{equation 8}), we procure the following corollary.

\begin{corollary}
For $n\in 
\mathbb{N}
$, we have%
\begin{equation*}
\int_{%
\mathbb{Z}
_{p}}u^{\eta }\left( 1-\eta \right) ^{n}d\mu _{-1}\left( \eta \right) =\frac{%
2}{u+1}+\frac{2}{u^{2}+u}+\frac{2}{u^{3}+u}H_{n}\left( -u^{-1}\right) .
\end{equation*}
\end{corollary}

\section{\textbf{Some identities on the Frobenius-Euler numbers}}

\qquad In this section, we develop Frobenius-Euler numbers, that is, we
derive some interesting and worthwhile relations for studying in Theory of
Analytic Numbers.

\qquad Now also, for $x\in \left[ 0,1\right] $, we rewrite definition of
Bernstein polynomials as follows:%
\begin{equation}
B_{k,n}\left( x\right) =\binom{n}{k}x^{k}\left( 1-x\right) ^{n-k}\text{,
where }n,k\in 
\mathbb{Z}
_{+}\text{.}  \label{equation 9}
\end{equation}

By expression of (\ref{equation 9}), we have the properties of symmetry of
Bernstein polynomials as follows:%
\begin{equation}
B_{k,n}\left( x\right) =B_{n-k,n}\left( 1-x\right) \text{, (for detail, see 
\cite{kim 19}).}  \label{equation 10}
\end{equation}

Thus, from Corollary 1, (\ref{equation 9}) and (\ref{equation 10}), we see
that%
\begin{eqnarray*}
\int_{%
\mathbb{Z}
_{p}}B_{k,n}\left( \eta \right) u^{\eta }d\mu _{-1}\left( \eta \right)
&=&\int_{%
\mathbb{Z}
_{p}}B_{n-k,n}\left( 1-\eta \right) u^{\eta }d\mu _{-1}\left( \eta \right) \\
&=&\binom{n}{k}\sum_{l=0}^{k}\binom{k}{l}\left( -1\right) ^{k+l}\int_{%
\mathbb{Z}
_{p}}u^{\eta }\left( 1-\eta \right) ^{n-l}d\mu _{-1}\left( \eta \right) \\
&=&\binom{n}{k}\sum_{l=0}^{k}\binom{k}{l}\left( -1\right) ^{k+l}\left( \frac{%
2}{u+1}+\frac{2}{u^{2}+u}+\frac{2}{u^{3}+u}H_{n-l}\left( -u^{-1}\right)
\right) \text{.}
\end{eqnarray*}

For $n$, $k\in 
\mathbb{Z}
_{+}$ with $n>k$, we compute%
\begin{eqnarray}
&&\int_{%
\mathbb{Z}
_{p}}B_{k,n}\left( \eta \right) u^{\eta }d\mu _{-1}\left( \eta \right)
\label{equation 11} \\
&=&\binom{n}{k}\sum_{l=0}^{k}\binom{k}{l}\left( -1\right) ^{k+l}\left( \frac{%
2}{u+1}+\frac{2}{u^{2}+u}+\frac{2}{u^{3}+u}H_{n-l}\left( -u^{-1}\right)
\right)  \notag \\
&=&\left\{ \QATOPD. . {\frac{2}{u+1}+\frac{2}{u^{2}+u}+\frac{2}{u^{3}+u}%
H_{n}\left( -u^{-1}\right) ,\text{ \ \ \ \ \ \ \ \ \ \ \ \ \ \ \ \ \ \ \ \ \
\ \ \ \ \ \ \ \ \ \ \ \ \ \ if }k=0,}{\binom{n}{k}\sum_{l=0}^{k}\binom{k}{l}%
\left( -1\right) ^{k+l}\left( \frac{2}{u+1}+\frac{2}{u^{2}+u}+\frac{2}{%
u^{3}+u}H_{n-l}\left( -u^{-1}\right) \right) ,\text{ if }k>0.}\right.  \notag
\end{eqnarray}

Let us take the fermionic $p$-adic $q$-integral on $%
\mathbb{Z}
_{p}$ on the Bernstein polynomials of degree $n$ as follows:%
\begin{eqnarray}
\int_{%
\mathbb{Z}
_{p}}B_{k,n}\left( \eta \right) u^{\eta }d\mu _{-1}\left( \eta \right) &=&%
\binom{n}{k}\int_{%
\mathbb{Z}
_{p}}\eta ^{k}\left( 1-\eta \right) ^{n-k}u^{\eta }d\mu _{-1}\left( \eta
\right)  \label{equation 12} \\
&=&\frac{2}{u+1}\binom{n}{k}\sum_{l=0}^{n-k}\binom{n-k}{l}\left( -1\right)
^{l}H_{l+k}\left( -u^{-1}\right) \text{.}  \notag
\end{eqnarray}

Consequently, by expression of (\ref{equation 11}) and (\ref{equation 12})$,$
we state the following Theorem:

\begin{theorem}
The following identity holds true:%
\begin{equation*}
\sum_{l=0}^{n-k}\binom{n-k}{l}\left( -1\right) ^{l}H_{l+k}\left(
-u^{-1}\right) =\left\{ \QATOPD. . {1+u^{-1}+u^{-2}H_{n}\left(
-u^{-1}\right) ,\text{ \ \ \ \ \ \ \ \ \ \ \ \ \ \ \ \ \ \ \ \ \ \ \ \ if }%
k=0\text{,}}{\sum_{l=0}^{k}\binom{k}{l}\left( -1\right) ^{k+l}\left(
1+u^{-1}+u^{-2}H_{n-l}\left( -u^{-1}\right) \right) ,\text{ if }k>0\text{.}%
}\right.
\end{equation*}
\end{theorem}

Let $n_{1},n_{2},k\in 
\mathbb{Z}
_{+}$ with $n_{1}+n_{2}>2k$. Then, we derive the followings%
\begin{eqnarray}
&&\int_{%
\mathbb{Z}
_{p}}B_{k,n_{1}}\left( \eta \right) B_{k,n_{2}}\left( \eta \right) u^{\eta
}d\mu _{-1}\left( \eta \right)  \notag \\
&=&\binom{n_{1}}{k}\binom{n_{2}}{k}\sum_{l=0}^{2k}\binom{2k}{l}\left(
-1\right) ^{2k+l}\int_{%
\mathbb{Z}
_{p}}\left( 1-\eta \right) ^{n_{1}+n_{2}-l}u^{\eta }d\mu _{-1}\left( \eta
\right)  \notag \\
&=&\left( \binom{n_{1}}{k}\binom{n_{2}}{k}\sum_{l=0}^{2k}\binom{2k}{l}\left(
-1\right) ^{2k+l}\left( \frac{2}{u+1}+\frac{2}{u^{2}+u}+\frac{2}{u^{3}+u}%
H_{n_{1}+n_{2}-l}\left( -u^{-1}\right) \right) \right)  \notag \\
&=&\left\{ \QATOPD. . {\frac{2}{u+1}+\frac{2}{u^{2}+u}+\frac{2}{u^{3}+u}%
H_{n_{1}+n_{2}}\left( -u^{-1}\right) \text{, \ \ \ \ \ \ \ \ \ \ \ \ \ \ \ \
\ \ \ \ \ \ \ \ \ \ \ \ \ \ \ \ \ \ \ \ \ \ \ \ \ \ \ \ if }k=0,}{\binom{%
n_{1}}{k}\binom{n_{2}}{k}\sum_{l=0}^{2k}\binom{2k}{l}\left( -1\right)
^{2k+l}\left( \frac{2}{u+1}+\frac{2}{u^{2}+u}+\frac{2}{u^{3}+u}%
H_{n_{1}+n_{2}-l}\left( -u^{-1}\right) \right) ,\text{ \ if }k\neq 0.}\right.
\notag
\end{eqnarray}

Therefore, we obtain the following Theorem:

\begin{theorem}
For $n_{1},n_{2},k\in 
\mathbb{Z}
_{+}$ with $n_{1}+n_{2}>2k,$ we have%
\begin{eqnarray*}
&&\int_{%
\mathbb{Z}
_{p}}B_{k,n_{1}}\left( \eta \right) B_{k,n_{2}}\left( \eta \right) u^{\eta
}d\mu _{-1}\left( \eta \right) \\
&=&\left\{ \QATOPD. . {\frac{2}{u+1}+\frac{2}{u^{2}+u}+\frac{2}{u^{3}+u}%
H_{n_{1}+n_{2}}\left( -u^{-1}\right) \text{, \ \ \ \ \ \ \ \ \ \ \ \ \ \ \ \
\ \ \ \ \ \ \ \ \ \ \ \ \ \ \ \ \ \ \ \ \ \ \ \ if }k=0,}{\binom{n_{1}}{k}%
\binom{n_{2}}{k}\sum_{l=0}^{2k}\binom{2k}{l}\left( -1\right) ^{2k+l}\left( 
\frac{2}{u+1}+\frac{2}{u^{2}+u}+\frac{2}{u^{3}+u}H_{n_{1}+n_{2}-l}\left(
-u^{-1}\right) \right) ,\text{ \ if }k\neq 0.}\right.
\end{eqnarray*}
\end{theorem}

By using the binomial theorem, we can derive the following equation.%
\begin{eqnarray}
&&\int_{%
\mathbb{Z}
_{p}}B_{k,n_{1}}\left( \eta \right) B_{k,n_{2}}\left( \eta \right) u^{\eta
}d\mu _{-1}\left( \eta \right)  \label{equation 14} \\
&=&\dprod\limits_{i=1}^{2}\binom{n_{i}}{k}\sum_{l=0}^{n_{1}+n_{2}-2k}\binom{%
n_{1}+n_{2}-2k}{l}\left( -1\right) ^{l}\int_{%
\mathbb{Z}
_{p}}\eta ^{2k+l}u^{\eta }d\mu _{-1}\left( \eta \right)  \notag \\
&=&\frac{2}{u+1}\dprod\limits_{i=1}^{2}\binom{n_{i}}{k}%
\sum_{l=0}^{n_{1}+n_{2}-2k}\binom{n_{1}+n_{2}-2k}{l}\left( -1\right)
^{l}H_{2k+l}\left( -u^{-1}\right) .  \notag
\end{eqnarray}

Thus, we can obtain the following Corollary:

\begin{corollary}
For $n_{1},n_{2},k\in 
\mathbb{Z}
_{+}$ with $n_{1}+n_{2}>2k,$ we have%
\begin{eqnarray*}
&&\sum_{l=0}^{n_{1}+n_{2}-2k}\binom{n_{1}+n_{2}-2k}{l}\left( -1\right)
^{l}H_{2k+l}\left( -u^{-1}\right) \\
&=&\left\{ \QATOPD. . {1+u^{-1}+u^{-2}H_{n_{1}+n_{2}}\left( -u^{-1}\right) 
\text{, \ \ \ \ \ \ \ \ \ \ \ \ \ \ \ \ \ \ \ \ \ \ \ \ \ \ \ \ \ if }%
k=0,}{\sum_{l=0}^{2k}\binom{2k}{l}\left( -1\right) ^{2k+l}\left(
1+u^{-1}+u^{-2}H_{n_{1}+n_{2}-l}\left( -u^{-1}\right) \right) ,\text{ \ if }%
k\neq 0\text{.}}\right.
\end{eqnarray*}
\end{corollary}

For $\eta \in 
\mathbb{Z}
_{p}$ and $s\in 
\mathbb{N}
$ with $s\geq 2,$ let $n_{1},n_{2},...,n_{s},k\in 
\mathbb{Z}
_{+}$ with $\sum_{l=1}^{s}n_{l}>sk$. Then we take the fermionic $p$-adic $q$%
-integral on $%
\mathbb{Z}
_{p}$ for the Bernstein polynomials of degree $n$ as follows:%
\begin{eqnarray*}
&&\int_{%
\mathbb{Z}
_{p}}\underset{s-times}{\underbrace{B_{k,n_{1}}\left( \eta \right)
B_{k,n_{2}}\left( \eta \right) ...B_{k,n_{s}}\left( \eta \right) }}u^{\eta
}d\mu _{-1}\left( \eta \right) \\
&=&\dprod\limits_{i=1}^{s}\binom{n_{i}}{k}\int_{%
\mathbb{Z}
_{p}}\eta ^{sk}\left( 1-\eta \right) ^{n_{1}+n_{2}+...+n_{s}-sk}u^{\eta
}d\mu _{-1}\left( \eta \right) \\
&=&\dprod\limits_{i=1}^{s}\binom{n_{i}}{k}\sum_{l=0}^{sk}\binom{sk}{l}\left(
-1\right) ^{l+sk}\int_{%
\mathbb{Z}
_{p}}\left( 1-\xi \right) ^{n_{1}+n_{2}+...+n_{s}-l}u^{\eta }d\mu
_{-1}\left( \eta \right) \\
&=&\left\{ \QATOPD. . {\frac{2}{u+1}+\frac{2}{u^{2}+u}+\frac{2}{u^{3}+u}%
H_{n_{1}+n_{2}+...+n_{s}}\left( -u^{-1}\right) ,\text{ \ \ \ \ \ \ \ \ \ \ \
\ \ \ \ \ \ \ \ \ \ \ \ \ \ \ \ \ \ \ \ \ \ \ \ \ \ \ \ \ \ \ if }%
k=0,}{\dprod\limits_{i=1}^{s}\binom{n_{i}}{k}\sum_{l=0}^{sk}\binom{sk}{l}%
\left( -1\right) ^{sk+l}\left( \frac{2}{u+1}+\frac{2}{u^{2}+u}+\frac{2}{%
u^{3}+u}H_{n_{1}+n_{2}+...+n_{s}-l}\left( -u^{-1}\right) \right) ,\text{ \ \
if }k\neq 0.}\right.
\end{eqnarray*}

So from above, we have the following Theorem:

\begin{theorem}
For $s\in 
\mathbb{N}
$ with $s\geq 2$, let $n_{1},n_{2},...,n_{s},k\in 
\mathbb{Z}
_{+}$ with $\sum_{l=1}^{s}n_{l}>sk$. Then we have%
\begin{eqnarray*}
&&\int_{%
\mathbb{Z}
_{p}}u^{\eta }\dprod\limits_{i=1}^{s}B_{k,n_{i}}\left( \eta \right) u^{\eta
}d\mu _{-1}\left( \eta \right) \\
&=&\left\{ \QATOPD. . {\frac{2}{u+1}+\frac{2}{u^{2}+u}+\frac{2}{u^{3}+u}%
H_{n_{1}+n_{2}+...+n_{s}}\left( -u^{-1}\right) ,\text{ \ \ \ \ \ \ \ \ \ \ \
\ \ \ \ \ \ \ \ \ \ \ \ \ \ \ \ \ \ \ \ \ \ \ \ \ \ \ \ \ if }%
k=0,}{\dprod\limits_{i=1}^{s}\binom{n_{i}}{k}\sum_{l=0}^{sk}\binom{sk}{l}%
\left( -1\right) ^{sk+l}\left( \frac{2}{u+1}+\frac{2}{u^{2}+u}+\frac{2}{%
u^{3}+u}H_{n_{1}+n_{2}+...+n_{s}-l}\left( -u^{-1}\right) \right) ,\text{ \ \
if }k\neq 0.}\right.
\end{eqnarray*}
\end{theorem}

From the definition of Bernstein polynomials and the binomial theorem, we
easily get%
\begin{eqnarray}
&&\int_{%
\mathbb{Z}
_{p}}\underset{s-times}{\underbrace{B_{k,n_{1}}\left( \eta \right)
B_{k,n_{2}}\left( \eta \right) ...B_{k,n_{s}}\left( \eta \right) }}u^{\eta
}d\mu _{-1}\left( \eta \right)  \notag \\
&=&\dprod\limits_{i=1}^{s}\binom{n_{i}}{k}\sum_{l=0}^{n_{1}+...+n_{s}-sk}%
\binom{\sum_{d=1}^{s}\left( n_{d}-k\right) }{l}\left( -1\right) ^{l}\int_{%
\mathbb{Z}
_{p}}\xi ^{sk+l}u^{\eta }d\mu _{-1}\left( \eta \right)  \notag \\
&=&\frac{2}{u+1}\dprod\limits_{i=1}^{s}\binom{n_{i}}{k}%
\sum_{l=0}^{n_{1}+...+n_{s}-sk}\binom{\sum_{d=1}^{s}\left( n_{d}-k\right) }{l%
}\left( -1\right) ^{l}H_{sk+l}\left( -u^{-1}\right) .  \label{equation 17}
\end{eqnarray}

Therefore, by (\ref{equation 17}), we get novel properties of
Frobenius-Euler numbers with the following corollary:

\begin{corollary}
For $s\in 
\mathbb{N}
$ with $s\geq 2$, let $n_{1},n_{2},...,n_{s},k\in 
\mathbb{Z}
_{+}$ with $\sum_{l=1}^{s}n_{l}>sk.$ Then, we have 
\begin{eqnarray*}
u^{2} &&\sum_{l=0}^{n_{1}+...+n_{s}-sk}\binom{\sum_{d=1}^{s}\left(
n_{d}-k\right) }{l}\left( -1\right) ^{l}H_{sk+l}\left( -u^{-1}\right) \\
&=&\left\{ \QATOPD. . {u^{2}+u+H_{n_{1}+n_{2}+...+n_{s}}\left(
-u^{-1}\right) ,\text{ \ \ \ \ \ \ \ \ \ \ \ \ \ \ \ \ \ \ \ \ \ \ \ \ \ \ \
\ \ if }k=0,}{\sum_{l=0}^{sk}\binom{sk}{l}\left( -1\right) ^{sk+l}\left(
u^{2}+u+H_{n_{1}+n_{2}+...+n_{s}-l}\left( -u^{-1}\right) \right) ,\text{ \ \
if }k\neq 0.}\right.
\end{eqnarray*}
\end{corollary}

\end{document}